\documentclass[a4paper,11pt]{article}

\usepackage{latexsym}
\usepackage{amsmath,amsfonts,amssymb,amsthm}
\newtheorem{thm}{Theorem}[section]
\newtheorem{lem}[thm]{Lemma}

\newtheorem{prop}[thm]{Proposition}

\newtheorem{rem}[thm]{Remark}

\input epsf

\title{An easy upper bound for Ramsey numbers}
\author{Roland Bacher}

\begin{document}
\maketitle

{\sl Abstract: We prove easy 
upper bounds for Ramsey numbers. \footnote{
Math. class: 5D10, 5C55}.
}

\section{Introduction}

A theorem of Ramsey, see \cite{R}, implies the existence of a
smallest natural integer $R(n)$, now called the \emph{$n-$th Ramsey number},
such that every (simple unoriented)
graph $G$ with at least $R(n)$ vertices
contains either a complete graph with $n$ vertices or $n$
pairwise non-adjacent vertices (defining a complete graph in the 
complementary graph of $G$). 

The aim of this paper is to give a new simple proof of 
the following upper bound for Ramsey numbers:

\begin{thm}\label{thmgraph} We have 
$$R(n)\leq 2^{2n-3}$$
for $n\geq 2$.
\end{thm}

The currently best asymptotic upper bound,
$$R(n+1)\leq {2n\choose n}n^{-C\log n/\log\log n}\ ,$$
(for a suitable constant $C$) is due to Conlon, see \cite{C}. 

The standard proof of Ramsey's theorem, due to Erd\"os and Szekeres 
(see \cite{ES} or Chapter 35 of \cite{AZ}),
uses a two parameter Ramsey number $R(a,b)$ defined as the smallest
integer such that every graph with $R(a,b)$ vertices contains either
a complete graph with $a$ vertices or a subset of $b$ non-adjacent 
vertices. It is slightly more involved than our proof and gives 
the upper bound $R(n+1)\leq {2n\choose n}$
(based on the trivial values $R(a,1)=R(1,a)=1$
and on the inequality $R(a,b)\leq R(a-1,b)+R(a,b-1)$ for 
$a,b>1$).

Simple graphs are equivalent to complete graphs with edges of two
colours (encoding edges, respectively nonedges of simple graphs).
There is a generalization of Ramsey's theorem to an arbitrary 
finite number $m$ of colours as follows:
There exists a smallest natural number $R_m(n)$
such that every complete graph on $R_m(n)$ vertices with 
edges of $m$ colours contains $n$ vertices belonging to
a complete edge-monochromatic subgraph. The following result gives an 
upper bound for $R_m(n)$:

\begin{thm}\label{thmcol} We have
$$R_m(n)\leq 1+\sum_{j=0}^{m(n-2)} m^j=1+\frac{m^{mn-2m+1}-1}{m-1}$$
for $m,n\geq 2$.
\end{thm}

For $m=2$, the upper bound $1+\frac{m^{mn-2m+1}-1}{m-1}$ of 
Theorem \ref{thmcol} coincides with the upper bound
$2^{2n-3}=1+\frac{2^{2n-4+1}-1}{2-1}=2^{2n-3}$ for $R_2(n)=R(n)$
given by Theorem \ref{thmgraph}.

This paper contains a simple proof of Theorem \ref{thmgraph} 
(Section \ref{sectionm2})
and Theorem \ref{thmcol} (Section \ref{sectcol}) which is a variation
on the proofs usually found and is perhaps slightly simpler.
In Section \ref{sectgen} we discuss a few
generalizations of the numbers $R'(n)$ and $R'_m(n)$ playing 
a crucial role in the proofs.

\section{Proof of Theorem \ref{thmgraph}}\label{sectionm2}

Given a finite graph $G$, we define $\rho'(G)$ to be the largest natural number
such that $G$ contains two 
(not necessarily disjoint) subsets $A$ and $B$ of vertices satisfying the 
following two conditions: \begin{enumerate}
\item{}
All vertices of $A$ are adjacent to each other and no vertices of $B$
are adjacent.
\item{} $\sharp(A)+\sharp(B)=\rho'(G)$.
\end{enumerate}

In this section, the letter $A$ always denotes a set of pairwise adjacent 
vertices and $B$ denotes a set of pairwise non-adjacent vertices.
Two such subsets $A,B$ of vertices in a graph $G$ realize $\rho'(G)$
if $\rho'(G)=\sharp(A)+\sharp(B)$.

We define $R'(n)$ as the smallest 
natural integer such that $\rho'(G)\geq n$ for every graph $G$ with $R'(n)$
vertices.

\begin{lem}\label{lemineqRRpr} We have $R(n)\leq R'(2n-1)$.
\end{lem}

\noindent{\bf Proof} A graph $G$ with $R'(2n-1)$ vertices
contains subsets $A$ and $B$ of vertices realizing
$\rho'(G)\geq 2n-1$. One of the subsets $A,B$ thus contains 
at least $n$ vertices. If $\sharp(A)\geq n$, 
the graph $G$ contains a complete subgraph of $n$ vertices,
if $\sharp(B)\geq n$, the graph $G$ contains 
$n$ pairwise non-adjacent vertices.\hfill$\Box$

\begin{lem}\label{lemineqRpr} We have $R'(n+1)\leq 2R'(n)$.
\end{lem}

\noindent{\bf Proof} We choose a vertex $v$
in a graph $G$ with $2R'(n)$ vertices. We denote by $G_v$
the subgraph of $G\setminus \{v\}$ defined by all neighbours
of $v$. Up to replacing $G$ by its 
complementary graph (and exchanging the roles of the sets
$A$ and $B$), we can suppose that $G_v$ has at least 
$\lceil (2R'(n)-1)/2\rceil=R'(n)$ vertices. 
Hence we have $\rho'(G)\geq n$ and we can find subsets 
$A,B$ of vertices in $G_v$ which realize $\rho'(G)$.
The subset $A\cup\{v\}$ contains thus $\sharp(A)+1$ pairwise adjacent
vertices of $G$ and we have 
$\rho'(G)\geq \sharp(A\cup\{v\})+\sharp(B)=\rho'(G_v)+1\geq n+1$.
\hfill$\Box$

\begin{prop}\label{propupperRp} We have $R'(n)\leq 2^{n-2}$ for $n\geq 2$.
\end{prop}

\noindent{\bf Proof} If $n=2$ we take $A=B=\{v\}$ 
where $v$ is the unique
vertex of the trivial graph $G=\{v\}$ on one vertex $v$.

Induction on $n$ using Lemma \ref{lemineqRpr} ends the proof.\hfill$\Box$

\noindent{\bf Proof of Theorem \ref{thmgraph}}
The proof follows from the inequalities
$$R(n)\leq R'(2n-1)\leq 2^{2n-1-2}=2^{2n-3}$$
given by Lemma \ref{lemineqRpr} and 
Proposition \ref{propupperRp}.\hfill$\Box$
 
\begin{rem} The proof of Proposition \ref{propupperRp} can easily 
be rewritten algorithmically: Given a graph $G$ with at least
$2^{n-2}\geq 2$ vertices,
set $A=B=\emptyset$. While $G$ has at least $4$ vertices, choose a vertex $v$.
Set $A=A\cup\{v\}$ and replace $G$ by the subgraph induced by all neighbours of $v$ if $v$ has more neighbours than non-adjacent vertices in $G$. Otherwise
set $B=B\cup\{v\}$ and replace $G$ by the subgraph induced by all 
non-neighbours ($\not=v$) of $v$.

If $G$ has $2$ or $3$ vertices, choose two vertices $v,w$ and replace $A$
by $A\cup\{v,w\}$ if $v$ and $w$ are adjacent. Otherwise replace $B$
by $B\cup\{v,w\}$.
\end{rem}

\subsection{Value of $R(3)$} The inequality $R'(4)\leq 2^{4-2}=4$
given by the case $n=2$ of Proposition \ref{propupperRp} is not sharp: Indeed,
we have $R'(4)=3$ as can be seen by inspecting all four possible graphs 
on three vertices. (The set $A$ has respectively $1,2,2,3$ elements for a 
$3-$vertex graph with $0,1,2,3$ edges.)
Lemma \ref{lemineqRpr} shows now $R'(5)\leq 6$ and we get
$R(3)\leq R'(5)\leq 6$ by Lemma \ref{lemineqRRpr}. Since a cycle
with $5$ vertices contains no triangle and no triplet of pairwise non-adjacent
vertices, both inequalities are sharp and we have $R(3)=R'(5)=6$.

\begin{rem} The value $R'(4)=3$ can of course be used for improving 
the upper bound $R(n)\leq 2^{2n-3}$ in Theorem \ref{thmgraph}
to $3\cdot 2^{2n-5}$ for $n\geq 3$. More generally, any interesting upper
bound on $R'(n)$ for $n>4$ easily yields an improvement of
Theorem \ref{thmgraph}.
\end{rem}

\section{Proof of Theorem \ref{thmcol}}\label{sectcol}

We define $R'_m(n)$ to be the smallest integer such that every
complete graph with $R'_m(n)$ vertices and
edges of $m$ colours contains $m$ (not necessarily disjoint)
subsets of vertices $A_1,\dots,A_m$ with $\sharp(A_1)+\sharp(A_2)+\dots+\sharp(A_m)=n$ and with $A_1,\dots,A_m$ defining $m$ complete edge-monochromatic 
graphs of different edge-colours.

Given a complete graph $G$ with $m-$coloured edges, we denote by $\rho'(G)$
the largest integer such that $G$ contains $m$ (not necessarily disjoint)
subsets $A_1,\dots,A_m$
of vertices defining complete edge-monochromatic subgraphs of 
different colours and $\rho'(G)=\sharp(A_1)+\dots+\sharp(A_m)$.
We say that $m$ such subsets $A_1,\dots,A_m$ realize $\rho'(G)$.

We have of course $\rho'(G)\geq n$ if $G$ contains at least $R'(n)$ vertices.

\noindent{\bf Examples:} \begin{enumerate}
\item{}
We have $R'_m(m)=1$ by setting $A_1=A_2=\dots=A_m=\{v\}$
where $v$ is the unique vertex of the trivial graph with one vertex (the 
empty sets of edges in $A_1,\dots,A_m$ have different colours by
convention).

The value $R'_m(m)=1$ also follows from $R_m(1)=1$ applied to the 
the trivial inequality $R'_m(n+m-1)\leq R_m(n)$ obtained by 
completing a complete edge-monochromatic subgraph on $n$ vertices
with $m-1$ singletons representing complete edge-monochromatic subgraphs
of the $m-1$ remaining colours.
\item{} $R'_m(m+1)=2$ since an edge-coloured complete graph on $2$
vertices is always monochromatic.
\item{} $R'_m(m+2)=3$ since every edge-coloured triangle is either 
edge-monochromatic or contains two edges of different colours.
\end{enumerate}

\begin{lem}\label{lemRmineq} We have $R_m(n)\leq R'(m(n-1)+1)$.
\end{lem}

\noindent{\bf Proof} A set of $m$ integers summing up to $m(n-1)+1$ 
contains an element at least equal to $n$. For every realization 
$A_1,\dots,A_m$ of $\rho'(G)\geq m(n-1)+1$ of a graph $G$ with 
$R'(m(n-1)+1)$ vertices there thus exists an index $i$ 
such that $A_i$ defines an edge-monochromatic complete 
graph on at least $n$ vertices.\hfill$\Box$

\begin{lem}\label{lemineqm} We have $R'_m(n+1)\leq 2+m(R'_m(n)-1)$.
\end{lem}

\noindent{\bf Proof} Fixing a vertex $v$ in a complete graph
$G$ with $2+m(R'_m(n)-1)$ vertices
and edges of $m$ colours, we get a partition $V\setminus \{v\}=
V_1\cup\dots\cup V_m$ of all vertices different from $v$ by
considering the set $V_i$ of vertices joined by an edge of colour $i$ to $v$.
Since $V\setminus\{v\}$ has $1+m(R'_m(n)-1)$ elements, there
exists a set $V_i$ containing at least $R'_m(n)$ vertices.
The subgraph $G_i$ with vertices $V_i$ thus contains a realization
$A_1,\dots,A_m$ of $\rho'_m(G_i)\geq n$. 
Since $v$ is joined by edges of colour $i$ to all elements of $A_i$ , 
the set of vertices $A_i\cup\{v\}$ defines a complete
edge-monochromatic subgraph of colour $i$ in $G$.
This proves $\rho'(G)\geq \sharp(A_1)+\dots+\sharp(A_i\cup\{v\})+\dots
+\sharp(A_m)=\rho'G_i)+1\geq n+1$.\hfill$\Box$

\begin{prop}\label{propupperRm} We have 
$$R'_m(m+k)\leq 1+\sum_{j=0}^{k-1}m^j= 1+\frac{m^k-1}{m-1}$$
for every natural integer $k$ (using the convention $\sum_{j=0}^{-1}m^j=0$ 
if $k=0$).
\end{prop}

\noindent{\bf Proof} The formula holds for $k=0$ with 
$A_1=A_2=\dots=A_m=\{v\}$ the unique vertex of the trivial graph $\{v\}$
reduced to one vertex.

Using Lemma \ref{lemineqm} and induction on $k$ we have
\begin{eqnarray*}
R'(m+k+1)&\leq& 2+m(R'_m(m+k)-1)\\
&\leq& 2+m\left(\left(1+\sum_{j=0}^{k-1}m^j\right)-1\right)\\
&=&1+\sum_{j=0}^km^j
\end{eqnarray*}
which ends the proof.\hfill$\Box$

\noindent{\bf Proof of Theorem \ref{thmcol}} We have
$$R_m(n)\leq R'(m(n-1)+1)\leq1+\frac{m^{m(n-1)+1-m}-1}{m-1}=1+
\frac{m^{mn-2m+1}-1}{m-1}$$
where the first inequality is Lemma \ref{lemRmineq}
and the second inequality is Proposition \ref{propupperRm}.
\hfill$\Box$

\section{Generalizations of the number $R'_m(n)$}\label{sectgen}

The number $R'_m(n)$ has two obvious generalizations.

The first one is given by considering
$R'_{m,j}(n)$ with $j\in \{1,\dots,m\}$ defined as the smallest integer
such that every complete graph with $R'_{m,j}(n)$
vertices and edges of $m$ colours
contains $j$ edge-monochromatic complete subgraphs 
of different edge-colours and of size $\alpha_1,\dots,\alpha_j$
such that $\alpha_1+\dots+\alpha_j=n$. Therefore we consider only
the $j$ colours corresponding to the $j$ largest edge-monochromatic
complete subgraphs. For $j=1$, we recover the usual Ramsey numbers 
$R_m(n)$, for $j=m$ we get the numbers $R'_m(n)$ introduced previously.

The second generalization depends on an unbounded function 
$s:\mathcal G\longrightarrow \mathbb N$ (one can also 
work with $m$ different unbounded
functions $s_c: \mathcal G\longrightarrow \mathbb N$
indexed by colours
or replace the target-set of natural integers by the set of non-negative 
real numbers) on the set $\mathcal G$ of all finite simple graphs.

For $n\geq 1$ we define $R'_{m,s}(n)$ 
as the smallest integer such that every complete graph on 
$R'_{m,s}(n)$ vertices contains $m$ (not necessarily complete)
edge-monochromatic subgraphs $G_1,\dots,G_m$ of colour
$1,\dots,m$ satisfying $s(G_1)+s(G_2)+\dots+s(G_m)\geq n$
(respectively $s_1(G_1)+\dots+s_m(G_m)\geq n$).

The numbers $R'_m(n)$ correspond to the choice 
$s(G)=n$ if $G$ is the complete graph on $n$ vertices and $s(G)=0$ otherwise.

Other perhaps interesting choices are $s(G)=n$ if $G$ is an $n$-cycle
and $s(G)=0$ otherwise, or $s(G)=n$ if $G$ is a simple path 
(two endpoints of degree $1$ and all other vertices of degree $2$)
with $n$ vertices.

It is of course possible to combine both generalizations by 
defining $R'_{m,j,s}(n)$ in the obvious way
considering only the $j$ colours giving the largest contributions
to the sum $s(G_1)+\dots+s(G_m)$.

\subsection{Analogues of $R'$ for van der Waerden numbers}
Van der Waerden's Theorem gives the existence of a function 
$W:\mathbb \{2,3,4,\dots\}\times \{2,3,4,\dots\}
\longrightarrow \mathbb N$ associating 
to two integers $m,n\geq 2$ the smallest natural integer
$W(m,n)$ such that every colouring of the $W(m,n)$ consecutive natural 
integers
$1,2,\dots,W(m,n)$ with $m$ colours contains a monochromatic arithmetic
progression with $n$ elements.

We define $W'(m,n)$ in the obvious way
as the smallest natural integer such that every 
colouring of $1,2,\dots,W'(m,n)$ with $m$ colours contains $m$
(perhaps empty) monochromatic progressions of different colours and 
of lengths
$\alpha_1,\dots,\alpha_m$ summing up to $n=\alpha_1+\dots+\alpha_m$.

We have $W(m,n)\leq W'(m,m(n-1)+1)$ since a set of $m$ integers
strictly smaller than $n$ sums up at most to $m(n-1)$.
It is easy to check that $W'(2,1)=1,W'(2,2)=2$ and $W'(2,3)=3$.

For $W'(2,4)$ we get $W'(2,4)=6$ as can be seen as follows:
$W'(2,4)>5$ by inspection of the black-white colouring $bbwbb$
of $1,2,3,4,5$. Consider a black-white colouring of $1,\dots,6$ 
not containing a black progression of size $\alpha$ and 
a white progression of size $\beta$ such that $\alpha+\beta\geq 4$.
Such a colouring cannot use only one colour (otherwise we 
can take $\alpha=6$ or $\beta=6$). It cannot use both colours twice
(otherwise we can take $\alpha=2$ and $\beta=2$).
It uses thus one colour, say white, only once and we have 
necessarily $\alpha\geq 3$ and $\beta=1$ since either $1,3,5$ 
(for an even white element) or $2,4,6$ (for an odd white element)
are all black.

It is of course also possible to consider the numbers $W'_j(m,n)$
defined by considering only the $j$ largest arithmetical progressions.
For $j=1$ we get the classical van der Waerden number $W(m,n)$.

\noindent{\bf Acknowledgements.} I thank S. Eliahou and G. McShane 
for comments.

\noindent Roland BACHER, Universit\'e Grenoble I, CNRS UMR 5582, Institut 
Fourier, 100 rue des maths, BP 74, F-38402 St. Martin d'H\`eres, France.

\noindent e-mail: Roland.Bacher@ujf-grenoble.fr

\end{document}